\documentclass{amsart}
\usepackage{amssymb}
\usepackage{ifthen}
\usepackage{graphicx}
\usepackage{setspace}

\newtheorem{thm}{Theorem}

\newtheorem{rem}{Remark}

\newcommand{\IC}{{\mathbb C}}
\newcommand{\ID}{{\mathbb D}}

\newcommand{\D}{{\mathbb D}}

\newcommand{\real}{{\operatorname{Re}}}

\def\be{\begin{equation}}
\def\ee{\end{equation}}
\newcommand{\bthm}{\begin{thm}}
\newcommand{\ethm}{\end{thm}}
\newcommand{\beqq}{\begin{eqnarray*}}
\newcommand{\eeqq}{\end{eqnarray*}}

\begin{document}
\title[On the initial coefficients for a class of analytic  functions]{On the initial coefficients for certain class of functions analytic in the unit disc}

\author[M. Obradovi\'{c}]{Milutin Obradovi\'{c}}
\address{Department of Mathematics,
Faculty of Civil Engineering, University of Belgrade,
Bulevar Kralja Aleksandra 73, 11000, Belgrade, Serbia}
\email{obrad@grf.bg.ac.rs}

\author[N. Tuneski]{Nikola Tuneski}
\address{Department of Mathematics and Informatics, Faculty of Mechanical Engineering, Ss. Cyril and Methodius
University in Skopje, Karpo\v{s} II b.b., 1000 Skopje, Republic of Macedonia.}
\email{nikola.tuneski@mf.edu.mk}


\subjclass[2000]{30C45, 30C50}
\keywords{analytic, univalent, coefficient bound, sharp results.}

\begin{abstract}
Let function $f$ be analytic in the unit disk $\ID$ and be normalized so that $f(z)=z+a_2z^2+a_3z^3+\cdots$. In this paper we give sharp bounds of the modulus of its second, third and fourth coefficient, if $f$ satisfies
\[ \left|\arg \left[\left(\frac{z}{f(z)}\right)^{1+\alpha}f'(z) \right] \right|<\gamma\frac{\pi}{2}   \quad\quad (z\in\ID),\]
for $0<\alpha<1$ and $0<\gamma\leq1$.
\end{abstract}


\maketitle

\section{Introduction and preliminaries}


Let ${\mathcal A}$ denote the family of all analytic functions
in the unit disk $\ID := \{ z\in \IC:\, |z| < 1 \}$ and  satisfying the normalization
$f(0)=0= f'(0)-1$.

A function $f\in \mathcal{A}$ is said to be {\it strongly starlike of order $\beta, 0<\beta \leq 1$ }if and only if
$$\left|\arg\frac{zf'(z)}{f(z)}\right|<\beta\frac{\pi}{2}\quad\quad (z\in\ID).$$
We denote this class by $\mathcal{S}^{\star}_{\beta}$. If $\beta=1$,
then $\mathcal{S}^{\star}_{1}\equiv \mathcal{S}^{\star}$ is the well-known class of \textit{starlike functions}.

In \cite{MO-1998} the author introduced the class $\mathcal{U}(\alpha, \lambda)$ ($0<\alpha$ and $\lambda <1$)
consisting of functions $f\in \mathcal{A}$  for which we have
$$
\left|\left(\frac{z}{f(z)}\right)^{1+\alpha}f'(z)-1\right|<\lambda \quad\quad (z\in\ID).
$$
In the same paper it is shown that $\mathcal{U}(\alpha, \lambda)\subset \mathcal{S}^{\star}$ if
$$0<\lambda \leq \frac{1-\alpha}{\sqrt{(1-\alpha)^{2}+\alpha^{2}}}.$$
The most valuable up to date results about this class can be found in Chapter 12 from \cite{DTV-book}.

In the paper \cite{MO-2000}  the author considered univalence of the class of functions $f\in \mathcal{A}$
satisfying the condition
\be\label{eq 17}
 \left|\arg \left[\left(\frac{z}{f(z)}\right)^{1+\alpha}f'(z) \right] \right|<\gamma\frac{\pi}{2}   \quad\quad (z\in\ID)
\ee
for $0<\alpha<1$ and $0<\gamma\leq1$, and proved the following theorem.

\medskip

\noindent {\it {\bf Theorem A}. Let $f\in \mathcal{A}$, $0<\alpha<2/\pi$ and let
$$ \left|\arg \left[\left(\frac{z}{f(z)}\right)^{1+\alpha}f'(z) \right]\right|<\gamma_{\star}(\alpha)\frac{\pi}{2} \quad\quad (z\in\ID),$$
where
$$ \gamma_{\star}(\alpha)=\frac{2}{\pi}\arctan \left(\sqrt{\frac{2}{\pi\alpha}-1}\right)-\alpha \sqrt{\frac{2}{\pi\alpha}-1}.$$
Then $f\in\mathcal{S}^{\star}_{\beta},$ where
$$\beta=\frac{2}{\pi}\arctan \sqrt{\frac{2}{\pi\alpha}-1}.$$ }

\section{Main result}

In this paper we will give the sharp estimates for initial  coefficients of functions $f\in \mathcal{A}$
which satisfied the condition \eqref{eq 17}. Namely, we have

\bthm\label{18-th 2}
Let $f(z)=z+a_{2}z^{2}+a_{3}z^{3}+\cdots$ belongs to the class $\mathcal{A}$ and satisfy the condition \eqref{eq 17} for $0<\alpha<1$ and $0<\gamma\leq1$. Then we have the  next sharp estimations:
  \smallskip
\begin{enumerate}
  \item[($a$)] $|a_{2}|\leq \frac{2\gamma}{1-\alpha}$;
  \smallskip
  \item[($b$)] $|a_{3}|\leq \left\{\begin{array}{cc}
\frac{2\gamma}{2-\alpha}, & 0<\gamma \leq \frac{(1-\alpha)^{2}}{3-\alpha}  \\[2mm]
\frac{2(3-\alpha)\gamma^{2}}{(1-\alpha)^{2}(2-\alpha)}, &\frac{(1-\alpha)^{2}}{3-\alpha} \leq \gamma \leq 1
\end{array}
\right.$;
\smallskip
  \item[($c$)] $|a_{4}|\leq \left\{\begin{array}{cc}
\frac{2\gamma}{3-\alpha},& 0< \gamma \leq \gamma_{\nu} \\[2mm]
\frac{2\gamma}{3-\alpha}\left[\frac{1}{3}+
\frac{2}{3}\frac{(\alpha ^{2}-6\alpha +17)\gamma^{2}}{(1-\alpha)^{3}(2-\alpha )}\right],&
\gamma_{\nu}\leq\gamma \leq1
\end{array}
\right . $;
\end{enumerate}
where
$$ \gamma_{\nu}=\sqrt{\frac{(1-\alpha)^{3}(2-\alpha )}{\alpha ^{2}-6\alpha +17}}.$$
\ethm

\begin{proof}
We can write the condition \eqref{eq 17}  in the form
\be\label{eq 18}
\left(\frac{f(z)}{z}\right)^{-(1+\alpha)}f'(z)=\left(\frac{1+\omega(z)}{1-\omega (z)}\right)^{\gamma}
\,\left (=(1+2\omega(z)+2\omega^{2}(z)+\cdots)^{\gamma}\right),
\ee
where $\omega$ is analytic in $\ID$ with $\omega(0)=0$ and $|\omega(z)|<1$, $z\in\ID$.
If we denote by $L$ and $R$ left and right hand side of equality \eqref{eq 18}, then we have
\beqq
\begin{split}
L&= \left[1-(1+\alpha)(a_{2}z+\cdots)+{-(1+\alpha)\choose 2}(a_{2}z+\cdots)^{2} \right.\\
&+\left.{-(1+\alpha)\choose3}(a_{2}z+\cdots)^{3}+\cdots\right] \cdot(1+2a_{2}z+3a_{3}z^{2}+4a_{4}z^{3}+\cdots)
\end{split}
\eeqq
and if we put $\omega(z)=c_{1}z+c_{2}z^{2}+\cdots$ :
\beqq
\begin{split}
R&= 1+\gamma\left[2(c_{1}z+c_{2}z^{2}+\cdots)+2(c_{1}z+c_{2}z^{2}+\cdots)^{2}+\cdots\right]\\
&+ {\gamma \choose 2}\left[2(c_{1}z+c_{2}z^{2}+\cdots)+2(c_{1}z+c_{2}z^{2}+\cdots)^{2}+\cdots\right]^{2}\\
&+ {\gamma \choose3} \left[2(c_{1}z+c_{2}z^{2}+\cdots)+2(c_{1}z+c_{2}z^{2}+\cdots)^{2}+\cdots\right]^{3}+\cdots.
\end{split}
\eeqq
If we compare the coefficients on $z,z^{2},z^{3}$ in $L$ and $R$, then, after some calculations, we obtain
\be\label{eq 19}
\begin{array}{l}
\displaystyle\smallskip
a_{2}=\frac{2\gamma}{1-\alpha}c_{1}, \\
\displaystyle\smallskip
a_{3}=\frac{2\gamma}{2-\alpha}c_{2}+\frac{2(3-\alpha)\gamma^{2}}{(1-\alpha)^{2}(2-\alpha)}c_{1}^{2},\\
\displaystyle\smallskip
a_{4}=\frac{2\gamma}{3-\alpha}\left(c_{3}+\mu c_{1}c_{2}+\nu c_{1}^{3}\right),
\end{array}
\ee
where
\be\label{eq 20}
\mu=\mu (\alpha,\gamma)=\frac{2(5-\alpha)\gamma}{(1-\alpha)(2-\alpha)}\quad\mbox{and}\quad
\nu=\nu (\alpha, \gamma)=\frac{1}{3}+
\frac{2}{3}\frac{(\alpha ^{2}-6\alpha +17)\gamma^{2}}{(1-\alpha)^{3}(2-\alpha )}.
\ee
Since $|c_{1}|\leq1 $, then by using \eqref{eq 19} we easily obtain the result (a) from this theorem.
Also , by using $|c_{1}|\leq1 $ and $|c_{2}|\leq1-|c_{1}|^{2} $, from \eqref{eq 19} we have
\beqq
\begin{split}
|a_{3}|&\leq\frac{2\gamma}{2-\alpha}|c_{2}|+\frac{2(3-\alpha)\gamma^{2}}{(1-\alpha)^{2}(2-\alpha)}|c_{1}|^{2}\\
&\leq\frac{2\gamma}{2-\alpha}(1-|c_{1}|^{2})+\frac{2(3-\alpha)\gamma^{2}}{(1-\alpha)^{2}(2-\alpha)}|c_{1}|^{2}\\
&= \frac{2\gamma}{2-\alpha}+\frac{2\gamma}{2-\alpha}\left[\frac{(3-\alpha)\gamma}{(1-\alpha)^{2}}-1\right]|c_{1}|^{2}
\end{split}
\eeqq
and the result depends of the sign of the factor in the last bracket.

The main tool of our proof for the coefficient $a_{4}$ will be the results of Lemma 2 in the paper  \cite{Prokhorov-1984}. Namely, in that paper the authors considered the sharp estimate of the functional $$\Psi(\omega)=|c_{3}+\mu c_{1}c_{2}+\nu c_{1}^{3}|$$ within the class of all holomorphic functions of the form
$$\omega(z)=c_{1}z+c_{2}z^{2}+\cdots $$
and satisfying the condition $|\omega(z)|<1$, $z\in\ID$. In the same paper in Lemma 2,
for $\omega$ of previous type and for any real numbers $\mu$ and $\nu$ they give the
sharp estimates $\Psi(\omega)\leq \Phi(\mu,\nu)$, where $\Phi(\mu,\nu)$ is given in general form in Lemma 2, and here we will use
$$\Phi(\mu,\nu) =  \left\{\begin{array}{cc}
1, & (\mu,\nu) \in D_1\cup D_2\cup \{(2,1)\}  \\[2mm]
|\nu|, & (\mu,\nu) \in \cup_{k=3}^7 D_k
\end{array}
\right.,$$
where
\smallskip

$ D_1 = \left\{(\mu,\nu):|\mu|\le\frac12,\, -1\le\nu\le1 \right\}, $
\smallskip

$ D_2 = \left\{(\mu,\nu):\frac12\le|\mu|\le2,\, \frac{4}{27}(|\mu|+1)^3-(|\mu|+1)\le\nu\le1 \right\}, $
\smallskip

$ D_3 = \left\{(\mu,\nu):|\mu|\le\frac12,\, \nu\le-1 \right\}, $
\smallskip

$ D_4 = \left\{(\mu,\nu):|\mu|\ge\frac12,\, \nu\le-\frac23(|\mu|+1) \right\}, $
\smallskip

$ D_5 = \left\{(\mu,\nu):|\mu|\le2,\, \nu\ge1 \right\}, $
\smallskip

$ D_6 = \left\{(\mu,\nu):2\le|\mu|\le4,\, \nu\ge\frac{1}{12}(\mu^2+8) \right\}, $
\smallskip

$ D_7 = \left\{(\mu,\nu):|\mu|\ge4,\, \nu\ge\frac23(|\mu|-1) \right\}. $
\smallskip

In that sense, we need the values $\alpha$ and $\gamma$ such that $0<\mu\leq \frac{1}{2}$, $\mu\leq2$,
$\mu\leq4$, $\nu\leq1.$ So, by using \eqref{eq 20}, we easily get the next equivalence
$$0<\mu\leq\frac{1}{2}\quad\Leftrightarrow \quad\gamma\leq \frac{(1-\alpha)(2-\alpha)}{4(5-\alpha)}:=\gamma_{1/2}\,;$$

$$\mu\leq 2\quad\Leftrightarrow \quad\gamma\leq \frac{(1-\alpha)(2-\alpha)}{5-\alpha}:=\gamma_{2}\,;$$

$$\mu\leq 4 \quad\Leftrightarrow \quad\gamma\leq \frac{2(1-\alpha)(2-\alpha)}{5-\alpha}:=\gamma_{4}\,;$$

$$ \nu \leq 1 \quad\Leftrightarrow \quad\gamma\leq\sqrt{\frac{(1-\alpha)^{3}(2-\alpha)}{\alpha ^{2}-6\alpha +17}}:=\gamma _{\nu}. $$
It is easily to obtain that all values $\gamma_{1/2}, \gamma_{2}, \gamma_{4}, \gamma_{\nu}$ are decreasing functions of $\alpha$, $0<\alpha<1$ and that
$$0<\gamma _{1/2}<\frac{1}{10},\quad0<\gamma _{2}<\frac{2}{5},\quad0<\gamma _{4}<\frac{4}{5},\quad
0<\gamma_{\nu}<\sqrt{\frac{2}{17}}=0.342997\ldots.$$
Also, it is clear that
$$0<\gamma_{1/2}<\gamma_{2}<\gamma_{4}$$
and it is easily to obtain that
$$\gamma_{1/2}\leq\gamma_{\nu} \quad\mbox{for}\quad\alpha\in (0, \alpha_{\nu}],$$
where $\alpha_{\nu}=0.951226\ldots$ is the root of the equation $5\alpha^{3}-56 \alpha^{2}+177 \alpha -122=0 $ (of course $\gamma_{\nu}\leq\gamma_{1/2}$ for $ \alpha\in [ \alpha_{\nu}, 1)$).

  \smallskip

Further, the next relation is valid:
$$0<\gamma_{\nu}<\gamma_{2}<\gamma_{4}.$$

\smallskip
\noindent{\bf Case 1.} ($0<\gamma\leq \gamma_{\nu}$). First, it means that $\nu \leq1$. If  $0<\gamma\leq \gamma_{1/2}$, then $0<\mu\leq\frac{1}{2}$ and $0<\nu\leq1$, which by Lemma 2 \cite{Prokhorov-1984} gives
$\Phi(\mu,\nu)=1$. If  $\gamma_{1/2}\leq \gamma\leq \gamma_{\nu}$, $\alpha \in(0,\alpha_{\nu})$, then $\frac{1}{2}\leq\mu<2$, $0<\nu\leq1$ and if we prove that
$$\frac{4}{27}(\mu+1)^{3}-(\mu+1)\leq \nu,$$
then also by Lemma 2 \cite{Prokhorov-1984} we have $\Phi(\mu,\nu)=1$. In that sense, let denote
$$L_{1}=: \frac{4}{27}(\mu+1)^{3}-(\mu+1)\quad\mbox{and}\quad R_{1}=\nu.$$
Since $L_{1}$ is an increasing function of $\mu$ for $\mu \geq\frac{1}{2}$ and since $\gamma \leq \gamma_{\nu}$, then
$$\mu \leq \frac{2(3-\alpha)\gamma_{\nu}}{(1-\alpha)(5-\alpha)}=2\sqrt{\frac{(1-\alpha)(5-\alpha)^{2}}{(2-\alpha)(\alpha ^{2}-6\alpha +17)}}<2\sqrt{\frac{25}{34}}=\frac{10}{\sqrt{34}}$$
(because the function under the square root is decreasing)
and so
$$L_{1}<\frac{4}{27}(\frac{10}{\sqrt{34}}+1)^{3}-(\frac{10}{\sqrt{34}}+1)=0.249838\ldots,$$
while
$$R_{1}=\nu =\frac{1}{3}+
\frac{2}{3}\frac{(\alpha ^{2}-6\alpha +17)\gamma^{2}}{(1-\alpha)^{3}(2-\alpha )}>\frac{1}{3}=0.33\ldots.$$
This implies the desired inequality.

\smallskip

\noindent{\bf Case 2.} ($\gamma_{\nu}\leq\gamma \leq 1 $) In this case we have that $\nu\geq1$. If
$\gamma_{\nu}\leq\gamma\leq \gamma_{2}$, $\alpha\in(0,1)$, then $0<\mu\leq2$, $\nu\geq 1$, which by Lemma 2 \cite{Prokhorov-1984} implies $\Phi(\mu,\nu)=\nu.$ If $\gamma_{2}\leq\gamma\leq \gamma_{4}$, $\alpha\in(0,1)$, then $2\leq\mu\leq4.$ Also, after some calculations, the inequality  $\nu\geq \frac{1}{12}(\mu^{2}+8)$ is equivalent to
$$\frac{43-23\alpha+5\alpha^{2}-\alpha^{3}}{(1-\alpha)^{3}(2-\alpha)^{2}}\gamma^{2}\geq1.$$
Since  $\gamma_{2}\leq\gamma$, then the previous inequality is satisfied if
$$\frac{43-23\alpha+5\alpha^{2}-\alpha^{3}}{(1-\alpha)^{3}(2-\alpha)^{2}}\gamma_{2}^{2}\geq1.$$
But, the last inequality is equivalent to the inequality $\alpha^{2}-2\alpha-3\leq0$, which is really true for $\alpha\in(0,1)$.
By Lemma 2 \cite{Prokhorov-1984} we also have $\Phi(\mu,\nu)=\nu .$  Finally, if $\gamma\geq \gamma _{4}$, then $\mu\geq4$ and if $\nu\geq \frac{2}{3}(\mu-1)$ we have (by using the same lemma) $\Phi(\mu,\nu)=\nu$. Really, the inequality $\nu\geq \frac{2}{3}(\mu-1)$ is equivalent with
$$2(\alpha ^{2}-6\alpha +17)\gamma^{2}-4(1-\alpha)^{2}(5-\alpha)\gamma +3(1-\alpha)^{3}(2-\alpha)\geq0.$$
Since the discriminant of previous trinomial is $$D=8(1-\alpha)^{3}(\alpha^{3}-2\alpha^{2}+17\alpha-52)<0$$
for $\alpha\in(0,1)$, then the previous inequality is valid.
By using \eqref{eq 19} we have that  $|a_{4}|\leq \frac{\gamma}{3-\alpha}$ (Case 1), or
$|a_{4}|\leq\frac{\gamma}{3-\alpha}\nu $ (Case 2), and from there the statement of the theorem.

All results of Theorem \ref{18-th 2} are the best possible as the functions $f_{i}$, $i=1,2,3$, defined with
$$\left(\frac{z}{f_{i}(z)}\right)^{1+\alpha}f_{i}'(z)=\left(\frac{1+z^i}{1-z^i}\right)^{\gamma},$$
where $0<\alpha<1$, $0<\gamma\leq1$. We have that
$$c_i=1 \quad\mbox{and}\quad c_j=0 \mbox{ when } j\neq i.$$
\end{proof}

\begin{rem}
By using Theorem A we can conclude that it is sufficient to be $\gamma \leq \gamma_{\star}(\alpha)$
and $0<\alpha <2/\pi$ for starlikeness  of functions $f\in \mathcal{A}$
which satisfied the condition \eqref{eq 17}.

Also, these conditions imply that the modulus of the coefficients $a_{2},$ $a_{3}$, $a_{4}$ is bounded with some constants. Namely, from the estimates given in Theorem \ref{18-th 2} we have, for example,
$$ |a_{2}|\leq \frac{2\gamma}{1-\alpha}\leq \frac{2\gamma_{\star}(\alpha)}{1-\alpha}, \quad\quad |a_{3}|\leq
\frac{2(3-\alpha)\gamma_{\star}^{2}(\alpha)}{(1-\alpha)^{2}(2-\alpha)},$$
etc.

We note that $\gamma_{\star}(\alpha)<1-\alpha $ for $0<\alpha <2/\pi$. Namely, if we put
$$\phi(\alpha)=:\gamma_{\star}(\alpha)-(1-\alpha ),$$
then $\phi'(\alpha)=1-\sqrt{\frac{2}{\pi\alpha}-1}$. It is easily to see that $\phi$ attains its minimum
$\phi(1/\pi)=-1/2$ and since $\phi(0+)=0, \phi(\frac{2}{\pi}-)=2/\pi-1<0$, we have the desired inequality.

When $\alpha\rightarrow 0$ , then $\gamma_{\star}(0+)=1$, and from Theorem \ref{18-th 2}, we have the next estimates for
$0<\gamma \leq1$:
\[
|a_{2}| \leq 2\gamma\leq2,\quad\quad |a_{3}| \leq \left\{\begin{array}{cc}
\gamma, & 0<\gamma \leq 1/3  \\[2mm]
  3\gamma^{2}, & 1/3 \leq \gamma \leq 1
\end{array}
\right.
\]
and
\[
|a_{4}| \leq \left\{\begin{array}{cc}
2\gamma/3, & 0<\gamma \leq \sqrt{2/17}  \\[2mm]
2\gamma(1+17\gamma^{2})/9\leq4 , & \sqrt{2/17}\le\gamma \leq 1
\end{array}
\right..
\]
This is the case when we have strongly starlike functions of order $\gamma$.

For $\gamma=1$ in Theorem \ref{18-th 2}, i.e. if $\real\left[\left(\frac{z}{f(z)}\right)^{1+\alpha}f'(z) \right]>0$, $z\in\D$, we have
\[|a_{2}|\leq \frac{2}{1-\alpha}, \quad\quad |a_{3}|\leq \frac{2(3-\alpha)}{(1-\alpha)^{2}(2-\alpha)}\]
and
\[|a_{4}|\leq \frac{2}{3-\alpha}\left[\frac{1}{3}+ \frac{2}{3}\frac{(\alpha ^{2}-6\alpha +17)}{(1-\alpha)^{3}(2-\alpha )}\right].\]
\end{rem}

\end{document}